\documentclass[journal,twoside,web]{ieeecolor}
\usepackage{generic}
\usepackage{cite}
\usepackage{amsmath,amssymb,amsfonts}
\usepackage{algorithmic}
\usepackage{graphicx}
\usepackage{textcomp}

\def\Re{\mathbb{R}}

\def\Lemma#1{Lemma~\ref{#1}}

\def\Sec#1{Sec.~\ref{#1}}

\def\notes#1{\marginpar{\tiny #1}\typeout{Notes!
Notes!
Notes!
}}
\renewcommand{\notes}[1]{\typeout{notes!}}
\def\FRAC#1#2#3{\genfrac{}{}{}{#1}{#2}{#3}}
\def\half{{\mathchoice{\FRAC{1}{1}{2}}%
{\FRAC{2}{1}{2}}%
{\FRAC{3}{1}{2}}%
{\FRAC{4}{1}{2}}}}

\newcommand{\tr}{\mbox{tr}}

\def\Re{\field{R}}

\def\Sec#1{Sec.~\ref{#1}}

\def\clP{{\cal P}}
\def\clZ{{\cal Z}}

\def\Sec#1{Sec~\ref{#1}}

\def\E{{\sf E}}

\def\Sec#1{Sec.~\ref{#1}}

\def\clZ{{\cal Z}}

\newtheorem{theorem}{Theorem}
\newtheorem{example}{Example}

\newtheorem{lemma}{Lemma}
\newtheorem{remark}{Remark}
\newtheorem{proposition}{Proposition}

\def\beq{\begin{eqnarray}} 
\def\bc{\begin{center}} 
\def\be{\begin{enumerate}}
\def\bi{\begin{itemize}} 
\def\bs{\begin{small}}
\def\bS{\begin{slide}}
\def\ec{\end{center}} 
\def\ee{\end{enumerate}}
\def\ei{\end{itemize}}
\def\es{\end{small}}
\def\eS{\end{slide}}
\def\eeq{\end{eqnarray}}

\newcommand{\ud}{\,\mathrm{d}}

\def\Re{\mathbb{R}}
\def\E{{\sf E}}

\def\clY{{\cal Y}}

\def\Sec#1{Sec.~\ref{#1}}
\def\Thm#1{Thm.~\ref{#1}}
\def\Prop#1{Prop.~\ref{#1}}

\def\clP{{\cal P}}
\def\clZ{{\cal Z}}

\renewcommand{\Re}{\mathbb{R}}

\def\FRAC#1#2#3{\genfrac{}{}{}{#1}{#2}{#3}}

\newcommand{\var}{\text{var}}

\def\clA{{\cal A}}

\def\clF{{\cal F}}

\def\clH{{\cal H}}

\def\clM{{\cal M}}

\def\clP{{\cal P}}

\def\clU{{\cal U}}
\def\clV{{\cal V}}

\def\clY{{\cal Y}}
\def\clZ{{\cal Z}}
\def\E{{\sf E}}
\def\bS{\mathbb{S}}
\def\sJ{{\sf J}}

\def\ones{{\sf 1}}
\def\sP{{\sf P}}
\def\tsP{{\tilde{\sf P}}}
\def\tE{{\tilde{\sf E}}}

\def\tp{{\hbox{\rm\tiny T}}}

\def\dv{\operatorname{diag}}

\def\dvar{\operatorname{var}}

\def\sQ{{\sf Q}}

\def\opt{{\text{\rm (opt)}}}

\def\divg{\nabla\cdot}

\begin{document}
\title{Duality for Nonlinear Filtering II: Optimal Control}
\author{Jin W. Kim, \IEEEmembership{Student member, IEEE}, and
  Prashant G. Mehta, \IEEEmembership{Senior member, IEEE}
	\thanks{This work is supported in part by the NSF award 1761622.}
	\thanks{Research reported in this paper was carried out by J. W. Kim, as part of his PhD
		dissertation work, while he was a graduate student at the University of Illinois at Urbana-Champaign.  
		He is now with the Institute of Mathematics at the University of Potsdam
		(e-mail: jin.won.kim@uni-potsdam.de).}
	\thanks{P. G. Mehta is with 
		the Coordinated Science Laboratory and the Department of Mechanical Science
		and Engineering at the University of Illinois at Urbana-Champaign
		(e-mail: mehtapg@illinois.edu).}}

\maketitle

\begin{abstract}
	This paper is concerned with the development and use of duality theory
	for a nonlinear filtering model with white noise observations.  The
	main contribution of this paper is to introduce a stochastic 
	optimal control problem as a dual to the nonlinear filtering problem.
	The mathematical statement of the dual relationship between the two
	problems is given in the form of a duality principle.  The
	constraint for the optimal control problem is the backward stochastic
	differential equation (BSDE) introduced in the companion 
	paper. The optimal control solution is obtained from an application of
	the maximum principle, and subsequently used to derive the equation of
	the nonlinear filter. The proposed duality is shown to be an exact
	extension of the classical Kalman-Bucy duality, and different
   from other types of optimal control and variational formulations given in literature.  
\end{abstract}


\section{Introduction}
\label{sec:introduction}

In this paper, we continue the development of duality theory for
nonlinear filtering.  While the companion paper (part I) was
concerned with a (dual) controllability counterpart of
stochastic observability, the purpose of this present paper (part
II) is to express the 
nonlinear filtering problem as a (dual) optimal control problem.  The proposed duality is shown to be an exact
extension of the original Kalman-Bucy duality~\cite{kalman1960general,kalman1961}, in the sense that the dual optimal control problem has the same minimum variance structure
for \emph{both} linear and nonlinear filtering problems.  Because of its
historical importance, we begin by introducing and 
reviewing the classical duality for the linear Gaussian model.

\subsection{Background and literature review}\label{ssec:lit-review}

The linear Gaussian filtering model is as follows:
\begin{subequations}\label{eq:linear-Gaussian-model}
	\begin{align}
		\ud X_t &= A^\tp X_t \ud t + \sigma \ud B_t,\quad X_0\sim N(m_0,\Sigma_0) \label{eq:linear-Gaussian-model-a}\\
		\ud Z_t &= H^\tp X_t \ud t + \ud W_t \label{eq:linear-Gaussian-model-b}
	\end{align}
\end{subequations}
where $X:=\{X_t\in\Re^d:0\leq t\leq T\}$
is the state process, the prior $N(m_0,\Sigma_0)$ is a Gaussian
density with mean $m_0\in \Re^d$ and variance $\Sigma_0 \succeq 0$,
$Z:=\{Z_t:0\leq t\leq T\}$ is the observation process, and both $B:=\{B_t:0\leq t\leq T\}$ and
$W:=\{W_t:0\leq t\leq T\}$ are Brownian motion (B.M.). It is
assumed that $X_0,B,W$ are mutually independent. 
The model parameters $A\in\Re^{d\times d}$, $H\in\Re^{d\times m}$, and $\sigma\in\Re^{d\times p}$.

For this problem, the dual optimal control formulations are
well-understood. These are of following two types: 

\begin{itemize}
\item Minimum variance optimal control problem:
\end{itemize}
\begin{subequations}\label{eq:mv-intro}
	\begin{align}
		\mathop{\text{Minimize}}_{\stackrel{\{u_t\in\Re^m:0\leq t\leq
			T\}}{=:u}}\!:\quad \sJ(u) &= |y_0|^2_{\Sigma_0}+ \int_0^{T} y_t^\tp
		(\sigma \sigma^\tp) y_t + |u_t|^2 \ud t \label{eq:mv-intro-a}\\
		\text{Subject to}\;\;:\; -\frac{\ud y_t}{\ud t}
		&= A y_t + H u_t,\quad y_T = f \;\; \text{(given)} \label{eq:mv-intro-b}
	\end{align}
\end{subequations}

\begin{itemize}
\item Minimum energy optimal control problem:
\end{itemize} 
\begin{subequations}\label{eq:mee-intro}
	\begin{align}
		\mathop{\text{Minimize}}_{\stackrel{\{u_t\in\Re^p:0\leq t\leq
				T\}=:u}{\tilde{m}_0\in\Re^d}}\!:\quad &
		\sJ(u,\tilde{m}_0;z)=
		|m_0-\tilde{m}_0|^2_{\Sigma_0^{-1}}
		\nonumber \\[-5pt]
		& \qquad + \int_0^{T} |u_t|^2+ |\dot{z}_t - H^\tp \tilde{m}_t|^2 \ud t \label{eq:mee-intro-a}\\
		\text{Subject to}\;\;:\;& \frac{\ud \tilde{m}_t}{\ud t}
		= A^\top \tilde{m}_t + \sigma u_t \label{eq:mee-intro-b}
	\end{align}
\end{subequations}
where $z = \{z_t\in\Re^m:0\leq t\leq T\}$ is a given sample path of observations.

\medskip

These two types of linear quadratic (LQ) optimal control problems are known since
1960s and described in~\cite[Sec.~7.3.1 and
7.3.2]{bensoussan2018estimation}.  Because it is discussed in the
seminal paper~\cite{kalman1961} of Kalman and Bucy, the minimum
variance duality~\eqref{eq:mv-intro} is also referred to as the
Kalman-Bucy duality~\cite{todorov2008general}.  The relationship of
the two problems to
the model~\eqref{eq:linear-Gaussian-model} is as follows:
\begin{itemize}
\item Minimum variance duality is related to the filtering problem for
  the model~\eqref{eq:linear-Gaussian-model}.
  The optimal control cost~\eqref{eq:mv-intro-a}
  comes from specifying a minimum variance objective
  for estimating the random variable $f^\tp X_T$ for $f\in\Re^d$.    
\item Minimum energy duality is related to a smoothing problem
  for the model~\eqref{eq:linear-Gaussian-model}. The optimal
  cost~\eqref{eq:mee-intro-a} is obtained from specifying a maximum
  likelihood (ML) objective for estimating a trajectory
  $\{\tilde{m}_t:0\leq t \leq T\}$ given a sample path $\{z_t:0\leq
  t\leq T\}$ of observations.
\end{itemize}
Their respective solutions are related to~\eqref{eq:linear-Gaussian-model} as follows:
\begin{itemize}
\item  The solution
of the minimum variance duality~\eqref{eq:mv-intro} is useful to
derive the Kalman filter
for~\eqref{eq:linear-Gaussian-model}~\cite[Ch.~7.6]{astrom1970}. The
derivation helps explain why the covariance equation of the Kalman
filter is the same as the differential Ricatti equation (DRE) of the
LQ optimal control.  
Note however that the time arrow is
reversed: the DRE is solved in forward time for the Kalman filter. 
This is because the constraint~\eqref{eq:mv-intro-b} is a
backward (in time) ordinary differential equation (ODE). 
\item The solution of the minimum
energy duality~\eqref{eq:mee-intro} is a favorite technique to derive
the forward-backward equations of smoothing for the
model~\eqref{eq:linear-Gaussian-model}.  The Hamilton's equation
for~\eqref{eq:mee-intro} is 
referred to as the Bryson-Frazier
formula~\cite[Eq.~(13.3.4)]{bryson2018applied}.  By introducing a DRE,
other forms of 
solution, e.g., the Fraser-Potter
smoother~\cite[Eq.~(16)-(17)]{fraser1969optimum}, are possible
and useful in practice.  
\end{itemize}

Given this background for the linear Gaussian model~\eqref{eq:linear-Gaussian-model}, there has been
extensive work spanning decades on extending duality to the problems of nonlinear
filtering and smoothing.  The prominent duality type solution
approaches in literature include the following:
\begin{itemize}
	\item Mortensen's maximum likelihood estimator (MLE)~\cite{mortensen1968}.
	\item Minimum energy estimator (MEE) 
        in the model predictive control (MPC)
	literature~\cite[Ch.~4]{rawlings2017model}.
	\item Log transformation relationship between the Zakai equation of
          nonlinear filtering and the Hamilton-Jacobi-Bellman (HJB)
          equation of optimal control~\cite{fleming1982optimal}.
	\item Mitter and Newton's variational formulation of the nonlinear
	smoothing problem~\cite{mitter2003}.
\end{itemize}
In an early work~\cite{mortensen1968}, Mortensen considered a slightly more general version
of the linear Gaussian model~\eqref{eq:linear-Gaussian-model} where the drift terms
in both~\eqref{eq:linear-Gaussian-model-a} and~\eqref{eq:linear-Gaussian-model-b} are
nonlinear.  Both the optimal control
problem and its forward-backward solution are straightforward
extensions of~\eqref{eq:mee-intro}.  Since 1960s, closely related
extensions have appeared by different names in different communities,
e.g., maximum likelihood estimation (MLE), maximum a posteriori (MAP)
estimation, and the minimum energy estimation (MEE) which is discussed next.


Based on the use of duality, the theory and algorithms developed in the MPC
literature are readily adapted to solve state estimation problems.  The
resulting class of estimators is referred to as the minimum energy
estimator (MEE)~\cite[Ch.~4]{rawlings2017model}.   
The MEE algorithms are broadly of two types: 
(i) Full information estimator (FIE) where the entire history of
observation is used; and (ii) Moving horizon estimator (MHE) where
only a most recent fixed window of observation is used.  An important
motivation is to also incorporate additional constraints in estimator
design.  
Early papers include~\cite{michalska1995moving,MHE01,krener2003convergence} and more recent
extensions have appeared in~\cite{Copp_Hespanha_simultaneous,farina_distributedMHE,mhe15_interacting,MHE_Switched}.
A historical survey is given 
in~\cite[Sec.~4.7]{rawlings2017model}  where Rawlings et. al. write
``{\em establishing duality [of optimal estimator] with the optimal
	regulator is a favorite technique for establishing estimator
	stability}''.  Although the specific comment is made for the
      Kalman filter, the remainder of the chapter amply demonstrates the
utility of dual constructions for {\em both} algorithm design and
convergence analysis (as the time-horizon $T\to\infty$).
Convergence analysis typically requires additional assumptions on
the model which in turn has motivated the work on nonlinear
observability and detectability definitions.  A literature survey of
these definitions, including the connections to duality theory,
appears in the introduction of the companion paper~\cite{duality_jrnl_paper_I}.    

While the focus of MEE is on deterministic models, duality is also an important theme in the study of
nonlinear stochastic systems (hidden Markov models).  
A key concept is the log
transformation~\cite{fleming1978exit}.
In~\cite{fleming1982optimal}, the log transformation
was used to transform the Zakai equation into
a Hamilton-Jacobi-Bellman (HJB) equation.
Because of this, the negative log of a posterior density is a
value function for some stochastic optimal control problem (this is
how duality is understood in stochastic settings~\cite[Sec.~4.8]{bensoussan1992stochastic}).
While the problem itself was not clarified
in~\cite{fleming1982optimal} (see
however~\cite{fleming1997deterministic}), Mitter and Newton introduced
a dual optimal control problem in~\cite{mitter2003} based on a
variational interpretation of the Bayes' formula.  This work continues to
impact algorithm design which remains an important area
of research~\cite{chen2016relation,kappen2016adaptive,reich2019data,ruiz_kappen2017,sutter2016variational}. 
A notable
ensuing contribution appeared in the PhD thesis-work~\cite{van2006filtering}
where Mitter-Newton duality is used to obtain results on nonlinear filter
stability.

Given the importance of duality for the purposes of stability analysis
in {\em both} deterministic and stochastic settings of the problem, it
is useful to return to the linear Gaussian
model~\eqref{eq:linear-Gaussian-model} and compare the two types of
duality~\eqref{eq:mv-intro} and~\eqref{eq:mee-intro}.  An important
point, that has perhaps not been stressed in literature, is that {\em minimum
  variance duality~\eqref{eq:mv-intro} is more compatible with the classical duality
  between controllability and observability in linear systems theory}.
This is because of the following reasons:
\begin{itemize}
\item{\em Inputs and outputs.} In~\eqref{eq:mv-intro}, the control input $u$ has the same
  dimension $m$ as the output process while in~\eqref{eq:mee-intro},
  the control input $u$ is the dimension $n$ of the process noise.  Evidently,
  it is natural to view the inputs and outputs as dual processes that
  have the same dimension.  
\item{\em Constraint.} If we ignore the noise terms
  in~\eqref{eq:linear-Gaussian-model} then the resulting
  deterministic state-output system
($\dot{x}_t=A^\top x_t$ and $z_t=H^\top x_t$)
shares a dual
  relationship with the deterministic state-input
  system~\eqref{eq:mv-intro-b}.  (It is shown in part
  I~\cite[Sec.~III-F]{duality_jrnl_paper_I}
  that~\eqref{eq:mv-intro-b} is also the dual for the stochastic
  system~\eqref{eq:linear-Gaussian-model}).
In contrast, the
  ODE~\eqref{eq:mee-intro-b} is a
  modified copy of the 
  model~\eqref{eq:linear-Gaussian-model-a}.  
\item{\em Stability condition.}  The condition for asymptotic
  analysis of~\eqref{eq:mv-intro} is stabilizability
  of~\eqref{eq:mv-intro-b} and by duality detectability of
  $(A^\tp,H^\tp)$.  
The latter is known to be also the appropriate condition 
  for stability of the Kalman filter. In contrast,
  for~\eqref{eq:mee-intro}, 
  asymptotic convergence of the optimal $\tilde{m}_T$ is possible even with $\sigma=0$.  The
  important 
  condition again is detectability of $(A^\tp,H^\tp)$ but it is not at
  all 
  easy to see from~\eqref{eq:mee-intro}.  
\item {\em Arrow of time.} Because the respective DREs are solved
  forward (resp. backward) in time for optimal filtering
  (resp. control), the
  arrow of time flips between optimal control and optimal
  filtering. Evidently, this is the case for minimum variance
  duality~\eqref{eq:mv-intro} but not so for the minimum energy
  duality~\eqref{eq:mee-intro}: The constraint~\eqref{eq:mv-intro-b} is
  a backward in time ODE while the constraint~\eqref{eq:mee-intro-b}
  is a modified copy of the signal model which proceeds forward in
  time.  
\end{itemize}  

All of this suggests that a fruitful approach -- for both defining
observability and for using the definition for asymptotic stability
analysis -- is to consider the minimum variance duality, which
naturally begets the following questions:
\begin{itemize}
\item What are the appropriate extensions of~\eqref{eq:mv-intro}
and~\eqref{eq:mee-intro} for nonlinear deterministic and stochastic systems?
\item What type of duality is implicit in Mitter-Newton's work?  It
is already evident that MEE is an extension of~\eqref{eq:mee-intro}.  
\end{itemize}
Both of these questions are answered in the present paper (for the
white noise observation model).  Before
discussing the original contributions, it is noted that the past work
on minimum variance duality has been on refinement and
extensions of the linear model with additional constraints.  
In~\cite{simon1970duality}, it is used to obtain the solution
to a class of singular regulator problems, and in~\cite{goodwin2005},
the Lagrangian dual for an MEE problem with truncated measurement
noise is considered.  Numerical algorithms for~\eqref{eq:mv-intro} and
its extensions appear
in~\cite{mishra2022automatica,mvFIR07,min_variance_fir_tv,mv_ui}. Prior
to our work, it was widely believed that the nonlinear extension of
minimum variance duality is not possible~\cite{todorov2008general}.

\subsection{Summary of original contributions}

The main contribution of this paper is to present a minimum variance
dual to the nonlinear filtering problem.  As in the companion paper
(part I), the nonlinear filtering problem is for the HMM with the
white noise observation model.  The mathematical statement of the dual
relationship between optimal filtering and optimal control is given in
the form of a duality principle (\Thm{thm:duality-principle}).  The
principle relates the value of the control problem to the variance of
the filtering problem.  The classical Kalman-Bucy
duality~\eqref{eq:mv-intro} is recovered as a special case for the
linear-Gaussian model~\eqref{eq:linear-Gaussian-model}. 

Two approaches are described to solve the optimal
control problem: (i) Based on the use of the stochastic maximum
principle to derive the Hamilton's equation (\Thm{thm:optimal-solution}); and 
(ii) Based on a martingale characterization (\Thm{thm:martingale}).  A
formula for the optimal control as a feedback control law is obtained
and used to derive the equation of the optimal nonlinear filter.  Our
duality is also related to Mitter-Newton duality with a side-by-side
comparison in Table~\ref{tb:comparison}.



This paper is drawn from the PhD thesis of the first
author~\cite{JinPhDthesis}.  A prior conference version appeared
in~\cite{kim2019duality}.  While the duality principle was already
stated in the conference paper, it relied on a certain
assumption~\cite[Assumption A1]{kim2019duality} which has now been proved.
Various formulae are stated more simply, e.g., the use of carr\'e du
champ operator to specify the running cost.  Issues related to
function spaces have been clarified to a large extent.  While the
conference version relied on the innovation process, the present 
version directly works with  the observation process.  Such a choice
is more natural for the problem at hand.  As a result, most of the
results and certainly their proofs are novel.  Comparison with
Mitter-Newton duality is also novel.

\subsection{Paper outline}

The outline of the remainder of this paper is as follows: 
The mathematical model and necessary background appears in
\Sec{sec:back}.  The dual
optimal control problem together with the duality principle and its relation
to the linear-Gaussian case is described in \Sec{sec:duality-principle}.
Its solution using the maximum principle and the martingale characterization appears in
\Sec{sec:standard-form} and \Sec{ssec:martingale},
respectively. Duality-based derivation of the equation of the
nonlinear filter appears in \Sec{sec:derivation-of-the-filter}.
A comparison with Mitter-Newton duality is contained in
\Sec{sec:comp}.  The paper closes with some conclusions and directions
for future work in
\Sec{sec:conc}.  All the proof are contained in the Appendix.

\section{Background}
\label{sec:back}

We briefly review the model and the notation as presented
in~\cite{duality_jrnl_paper_I}.  Although the
presentation is self-contained, it is in an abbreviated form with a
focus on additional new concepts that are necessary for this paper.  

On the probability space $(\Omega, \clF_T, \sP)$, we consider a pair of continuous-time stochastic processes $(X,Z)$
as follows:
\begin{itemize}
	\item The \emph{state process} $X = \{X_t:\Omega\to \bS:0\le t \le
          T\}$ is a Feller-Markov process taking values in the state-space
          $\bS$. The prior is denoted by $\mu \in
          \clP(\bS)$ (space of probability measures) and $X_0\sim \mu$. The infinitesimal generator is
          denoted by $\clA$.	
\item  The \emph{observation process} $Z = \{Z_t:0\le t \le T\}$ satisfies the stochastic differential
          equation (SDE):
	\begin{equation}\label{eq:obs-model}
		Z_t = \int_0^t h(X_s) \ud s + W_t,\quad t \ge 0
	\end{equation}
	where $h:\bS\to \Re^m$ is referred to as the observation function and $W =
        \{W_t:0\le t \le T\}$ is an $m$-dimensional Brownian motion
        (B.M.). We write $W$ is $\sP$-B.M. 
	It is assumed that $W$ is independent of $X$.
\end{itemize}
The above is referred to as the \emph{white noise observation model} of 
nonlinear filtering. The model is denoted by 
$(\clA,h)$.

An important additional concept in this paper is the \emph{carr\'e 
	du champ} operator $\Gamma$ defined as follows (see~\cite{bakry2013analysis}):
\[
	(\Gamma f)(x) = (\clA f^2)(x) - 2f(x)(\clA f)(x),\quad x\in \bS
\]
where $f:\bS\to\Re$ is a test function.  Explicit
formulae for the most important examples are described next.

\subsection{Guiding examples}

\begin{example}[Finite state-space]
\label{ex:finite}
$\bS=\{1,2,\ldots,d\}$.  A real-valued function $f$ is
  identified with a vector in $\Re^d$ where the $i^{\text{th}}$ element of the vector is $f(i)$. In this manner, the generator
  $\clA$ of the Markov process is identified with a row-stochastic
  rate matrix $A\in\Re^{d\times d}$ (the non-diagonal elements of $A$
  are non-negative and the row sum is zero). 
The carr\'e du champ operator $\Gamma:\Re^d\to \Re^d$ is as follows:
\begin{equation}\label{eq:Gamma-finite}
	(\Gamma f)(i) = \sum_{j \in \mathbb{S}} A(i,j) (f(i) - f(j))^2,\quad i\in\bS
\end{equation}
\end{example}

\medskip

\begin{example}[Euclidean state-space] 
\label{ex:itodiffusion}
$\bS=\Re^d$.  
The Markov process $X$ is an It\^o diffusion modeled using a
stochastic differential equation (SDE):
\begin{equation*}\label{eq:dyn_sde}
	\ud X_t = a (X_t) \ud t + \sigma(X_t) \ud B_t,\quad X_0\sim\mu
\end{equation*}
where 
$a\in C^1(\Re^d;
\Re^d)$ and $\sigma\in C^2(\Re^d; \Re^{d\times p})$ satisfy
appropriate technical conditions such that a strong solution exists
for $[0,T]$, and $B=\{B_t:0\le t \le T\}$ is a standard B.M.~assumed
to be independent of $X_0$ and $W$.  In the Euclidean case, all the
measures are identified with their density.  In particular, we use the
notation $\mu$ to denote the probability density function of the prior.   

The infinitesimal generator $\clA$ acts on $C^2(\Re^d;\Re)$ functions in its 
domain according to~\cite[Thm. 7.3.3]{oksendal2003stochastic}
\[
(\clA f)(x):= a^\tp(x) \nabla f(x) + \half  \tr\big(\sigma\sigma^\tp(x)(D^2f)(x)\big),\quad x\in\Re^d
\]
where $\nabla f$ is the gradient vector and $D^2 f$ is
the Hessian matrix. For $f\in C^1(\Re^d;\Re)$, the carr\'e du champ operator is given by 
\begin{equation}\label{eq:Gamma-Euclidean}
	(\Gamma f) (x) = \big|\sigma^\tp(x) \nabla f(x) \big|^2,\quad x\in\Re^d
\end{equation}
\end{example}

\medskip

\begin{example}[Linear Gaussian model] The 
  model~\eqref{eq:linear-Gaussian-model} introduced
  in~\Sec{sec:introduction} is a special case of  It\^o diffusion
  where the drift terms are linear $a(x)=A^\tp x$ and $h(x)=H^\tp x$,
  the coefficient of the process noise $\sigma(x)=\sigma$ is a
  constant matrix, and the prior $\mu$ is a Gaussian density. 
A real-valued linear function is expressed as \[f(x) = 
\tilde{f}^\tp x,\quad x\in\Re^d\] 
where $\tilde{f}\in \Re^d$. Then $\clA f$ is also a linear function 
given by
\[
\big(\clA f\big)(x) = (A\tilde{f})^\tp x,\quad x\in\Re^d
\]
and $\Gamma f$ is a constant function given by
\begin{equation}\label{eq:Gamma-LG}
\big(\Gamma f\big)(x) = \tilde{f}^\tp \big(\sigma \sigma^\tp\big) \tilde{f},\quad x\in\Re^d 
\end{equation}
\end{example}

\subsection{Background on nonlinear filtering}

The canonical filtration $\clF_t = \sigma\big(\{(X_s,W_s):0\le s \le
t\}\big)$. The filtration generated by the observation is denoted by
$\clZ :=\{\clZ_t:0\le t\le T\}$  where $\clZ_t = \sigma\big(\{Z_s:0\le
s\le t\}\big)$.  A standard approach 
is based upon the Girsanov change of measure.
Suppose the model satisfies the Novikov's condition: $
\E\left(\exp\big(\half \int_0^T |h(X_t)|^2\ud t\big)\right) < \infty
$. 
Define a new measure $\tsP$
on $(\Omega,\clF_T)$ as follows:
\[
\frac{\ud \tsP}{\ud \sP} = \exp\Big(-\int_0^T 
h^\tp(X_t) \ud W_t - \half \int_0^T |h(X_t)|^2\ud t\Big) =: D_T^{-1}
\]
Then it is shown that the probability law for $X$ is unchanged but
$Z$ is a $\tsP$-B.M.~that is independent of $X$~\cite[Lem.~1.1.5]{van2006filtering}. The
expectation with respect to $\tsP$ is denoted by $\tE(\cdot)$.  

The two probability measures are used to define the un-normalized and
the normalized (or nonlinear) filter are as follows:  For $0\le t \le T$ and $f\in
C_b(\bS)$,
\begin{align*}
\text{(un-normalized filter)} \quad \sigma_t(f) & := \tE\big(D_tf(X_t)|\clZ_t\big)\\
\text{(nonlinear filter)} \quad \pi_t(f) & := \E\big(f(X_t)|\clZ_t\big)
\end{align*}
As the name suggests, $\pi_t(f) = \frac{\sigma_t(f)}{\sigma_t(\ones)}$ 
which is referred to as the 
Kallianpur-Striebel formula~\cite[Thm.~5.3]{xiong2008introduction}
(here $\ones$ is the constant function $\ones(x)=1$ for all $x\in\bS$). Combining the tower property of conditional expectation with the
change of measure gives
\begin{equation}\label{eq:COM_formula}
\E (f(X_t)) = \E (\pi_t(f)) =\tE (\sigma_t(f))
\end{equation}

\subsection{Function spaces}

The notation $L^2_{\clZ_T}(\Omega;\Re^m)$ and
$L^2_{\clZ}\big([0,T];\Re^m\big)$ is used to denote the Hilbert space
of $\clZ_T$-measurable random
vector and $\clZ$-adapted stochastic process, respectively.  
These Hilbert spaces suffice if the state-space is finite.   In general
settings, let $\clY$ denote a suitable Banach space of real-valued
functions on $\bS$, equipped with the norm $\|\cdot\|_\clY$. Then
\begin{itemize}
\item For a random function, the
Banach space $L^2_{\clZ_T}(\Omega;\clY) := \big\{F:\Omega\to
\clY: F\text{ is }\clZ_T\text{-measurable},\;
\tE\big(\|F\|_\clY^2\big) < \infty\big\}$.  
\item 
For a function-valued
stochastic process, the Banach space is
$L^2_{\clZ}([0,T];\clY) := \Big\{Y:\Omega\times
[0,T]\to \clY\;:\; Y\text{ is }\clZ\text{-adapted},\;
\tE\Big(\int_0^T\|Y_t\|_\clY^2 \ud t\Big) < \infty\Big\}$. 
\end{itemize}

In the remainder of this paper, we set $\clY := C_b(\bS)$ (the space of
continuous and bounded functions) equipped with
the sup-norm. The dual space $\clM(\bS)$ (the space of rba measures) is denoted by
$\clY^\dagger$ where the duality pairing $\langle f, \rho \rangle =
\rho(f)$ for $f\in\clY$ and $\rho\in \clY^\dagger$.   

\section{Main result: The duality principle}\label{sec:duality-principle}

\subsection{Problem statement}


For a function $F\in L^2_{\clZ_T}\big(\Omega;\clY\big)$, the
nonlinear filter $\pi_T(F)$ is the minimum variance estimate of $F(X_T)$~\cite[Sec.~6.1.2]{bensoussan2018estimation}:
\begin{equation}\label{eq:minimum-variance}
\pi_T(F) = \mathop{\operatorname{argmin}}_{S_T\in L^2_{\clZ_T}(\Omega;\Re)} \E\big(|F(X_T)-S_T|^2\big)
\end{equation}
Our goal is to express the above minimum variance optimization problem as a dual optimal control problem.

The conditional variance
is denoted by 
\begin{equation*}
\clV_T(F)  
 := \E\big(|F(X_T)-\pi_T(F)|^2|\clZ_T\big) 
            = \pi_T(F^2) - \big(\pi_T(F)\big)^2
\end{equation*}
For notational ease, the expected value of the conditional
variance is denoted by 
\[
\dvar_T(F) := \E\big(\clV_T(F)\big)
\]
Strictly speaking, the above is variance only at time $T=0$.  However,
the verbiage is consistent with the ``minimum variance''
interpretation of the nonlinear filter.

\subsection{Dual optimal control problem} 

The function space of admissible control inputs is denoted by $\clU :=
L^2_{\clZ}\big([0,T];\Re^m\big)$.
An element of $\clU$ is denoted $U = \{U_t: 0\le t\le T\}$.  
It is referred to as the control input.
The main contribution of this paper is the following problem. 

\begin{itemize}
\item Minimum variance optimal control problem:
\end{itemize}
\begin{subequations}\label{eq:dual-optimal-control}
	\begin{align}
&           \mathop{\text{Minimize:}}_{U\in\;\clU}\; \sJ_T(U)  =
\dvar_0(Y_0) + 
\E \Big(\int_0^T l (Y_t,V_t,U_t\,;X_t) \ud t \Big)\label{eq:dual-optimal-control-a}\\
&          \text{Subject to (BSDE constraint):} \nonumber \\ 
&  -\!\ud Y_t(x) = \big((\clA Y_t)(x) + h (x) (U_t +
            V_t(x))\big)\ud t - V_t^\tp(x)\ud Z_t \nonumber \\
&    \quad \;  Y_T (x)  = F(x), \;\; x \in \bS 
\label{eq:dual-optimal-control-b}
	\end{align}
\end{subequations}
where the running cost 
\[
l(y,v,u;x):= (\Gamma y)(x) + |u+v(x)|^2
\] 
and $\dvar_0(Y_0)=\E(|Y_0(X_0) - \mu(Y_0)|^2) = \mu(Y_0^2) -\mu(Y_0)^2$.

\medskip

\begin{remark}
The BSDE~\eqref{eq:dual-optimal-control-b} is introduced in the
companion paper (part I) as the dual control system.  The data for the BSDE is the
given terminal condition $F\in L^2_{\clZ_T}\big(\Omega;\clY\big)$ and
the control input $U\in \clU$.  The solution of the BSDE is the pair 
$(Y,V)=\{(Y_t,V_t):0\leq t\leq T\} \in
L^2_{\clZ}\big([0,T];\clY\times \clY^m\big)$ which is
(forward) adapted to the filtration $\clZ$.
Existence, uniqueness, and regularity theory for linear BSDEs is
standard and 
throughout the paper, we assume that the solution of BSDE $(Y,V)$ is
uniquely determined in
$L^2_{\clZ}\big([0,T];\clY\times \clY^m\big)$ for
each given $Y_T\in L^2_{\clZ_T}(\Omega;\clY)$ and $U\in
L_\clZ^2\big([0,T];\Re^m\big)$.  The well-posedness results for finite
state-space can be found in~\cite[Ch.~7]{yong1999stochastic} and for
the Euclidean state space 
in~\cite{ma1999linear}. 
\end{remark}


\medskip

The relationship of~\eqref{eq:dual-optimal-control} to the minimum
variance objective~\eqref{eq:minimum-variance} is given the following theorem. 

\medskip

\begin{theorem}[Duality principle]\label{thm:duality-principle}
	For any admissible control $U\in \clU$, consider an estimator
	\begin{equation}\label{eq:estimator}
		S_T := \mu(Y_0) - \int_0^T U_t^\tp \ud Z_t
	\end{equation}
	Then 
	\begin{equation}\label{eq:duality-principle}
		\sJ_T(U) = \E\big(|F(X_T)-S_T|^2\big)
	\end{equation}
	
\end{theorem}

\medskip

\begin{proof}
See Appendix~\ref{ss:pf-thm51}.
\end{proof}

\medskip


The problem~\eqref{eq:dual-optimal-control} is a
stochastic linear quadratic optimal control problem for which there is a well
established existence-uniqueness theory for the optimal control
solution.  Application of this theory is the subject of the following
section.  For now, we assume that the optimal control is well-defined
and denote it as $U^\opt = \{U_t^\opt:0\le t \le T\}$.
Because the right-hand side of the
identity~\eqref{eq:duality-principle} is bounded below by $\var_T(F)$, the duality gap
\[
	\sJ_T(U^\opt) - \dvar_T(F) \ge 0
\]
In order to conclude that the duality gap is zero, it is both
necessary and sufficient to show that there exists a $U\in \clU$ such
that the estimator $S_T$, as given by~\eqref{eq:estimator}, equals $\pi_T(F)$. Since $Z$ is a $\tsP$-B.M., the following lemma is a
consequence of the It\^o representation theorem for Brownian motion~\cite[Thm.~4.3.3]{oksendal2003stochastic}.

\medskip

\begin{lemma}\label{lem:representation}
	For any $F \in L_{\clZ_T}^2(\Omega;\clY)$, there exists a
        unique $U\in \clU$ such that
	\[
	\pi_T(F) = \tE\big(\pi_T(F)\big) - \int_0^T U_t^\tp \ud Z_t,\quad \tsP\text{-a.s.}
	\]
\end{lemma}

\medskip

\begin{proof}
See Appendix~\ref{ss:pf-lem51}.
\end{proof}

\medskip

Because the duality gap is zero, the following implications are to be had:
	\begin{itemize}
		\item The optimal control $U^\opt$ 
                  gives the
                  conditional mean
		\[
		\pi_T(F) = \mu(Y_0) - \int_0^T \big(U_t^\opt\big)^\tp \ud Z_t,\quad \sP\text{-a.s.}
		\]
		\item The optimal value is the expected value of the conditional variance
		\[
		\dvar_T(F) = \dvar_0(Y_0) + \E \Big(\int_0^T
                l(Y_t,V_t,U_t^\opt;X_t) \ud t\Big)
		\]
	where $(Y,V)$ is the optimally controlled stochastic process obtained with $U=U^\opt$
        in~\eqref{eq:dual-optimal-control-b}.
	\end{itemize}

In fact, these two implications carry over to the entire optimal trajectory. 

\medskip

\begin{proposition}\label{prop:optimal-solution}
	Suppose
	$U^\opt$ is the optimal control input and that $(Y,V)$ is the associated solution of the BSDE~\eqref{eq:dual-optimal-control-b}.  Then for almost every $0\le t \le T$,
	\begin{align}
		\pi_t(Y_t) &= \mu(Y_0) - \int_0^t \big(U_s^\opt \big)^\tp\ud Z_s,\quad \sP\text{-a.s.} \label{eq:estimator-t} \\
		\dvar_t(Y_t) &= \dvar_0(Y_0) +\E\Big(\int_0^t l(Y_s,V_s,U_s^\opt;X_s) \ud s\Big) \label{eq:estimator-t-variance} 			
	\end{align}
\end{proposition}

\medskip

\begin{proof}
See Appendix~\ref{ssec:pf-optimal-solution}.
\end{proof}
\medskip

Consequently, the expected value of the conditional variance is the
optimal cost-to-go (for a.e. $0\leq t\leq T$).  We do not yet have a
formula for the optimal control $U^\opt$. 
The difficulty arises because there is no HJB equation for
BSDE-constrained optimal control problem. Instead, the literature on
such problem utilizes the stochastic maximum principle for BSDE which 
is the subject of the next section.  Before that, we discuss the
linear Gaussian case.



	
	



\subsection{Linear Gaussian case}
\label{sec:LG-case}

The goal is to show that the classical Kalman-Bucy
duality~\eqref{eq:mv-intro} described in~\Sec{sec:introduction} for
the linear Gaussian model~\eqref{eq:linear-Gaussian-model} is a
special case.  
Consider a linear function $F(x) = f^\tp x$ where $f\in\Re^d$ is a
given deterministic vector.  The problem is to compute a minimum
variance estimate of the scalar random variable $f^\tp X_T$. It is given by
$\E(f^\tp  X_T|\clZ_T)$.  Now, it is a standard result in the theory
of Gaussian
processes that conditional expectation can be evaluated in the form of
a linear predictor~\cite[Cor.~1.10]{le2016brownian}.  For this reason,
it suffices to consider an estimator of the form 
\[
	S_T := b - \int_0^T u_t^\tp \ud Z_t
\]
where $b\in\Re$ and $u = \{u_t \in \Re^m: 0\le t\le T\}$ are both
deterministic (the lower case notation is used to stress this).  Consequently, for linear Gaussian estimation, we can restrict the admissible
space of control inputs to $L^2\big([0,T];\Re^m\big)$ which is a much
smaller subspace of $L_{\clZ}^2\big([0,T];\Re^m\big)$.  
Using a deterministic control $u$, and the terminal condition $F(x) = f^\tp x$, the solution of the BSDE is given by
\[
Y_t(x) = y_t^\tp x, \quad V_t(x)=0, \quad x\in \Re^d, \;\;0\leq t\leq T
\] 
where $y=\{y_t\in\Re^d:0\leq t\leq T\}$ is a solution of the backward
ODE:
		\begin{equation*}\label{eq:LTI-ctrl}
		-\frac{\ud y_t}{\ud t}
		= A y_t + H u_t,\quad y_T = f
		\end{equation*}
Using the formula~\eqref{eq:Gamma-LG} for the carr\'e du champ,
the running cost
\begin{align*}
l (Y_t, V_t, U_t; X_t) &= (\Gamma Y_t)(X_t) + |U_t+V_t(X_t)|^2 \\
& = y_t^\tp (\sigma\sigma^\tp) y_t + |u_t|^2
\end{align*} 
With the
Gaussian prior, the initial 
cost $\var_0(y_0) = y_0^\tp \Sigma_0 y_0$.
Combining all of the above, the optimal control problem~\eqref{eq:dual-optimal-control}
reduces to~\eqref{eq:mv-intro} for the linear Gaussian
model~\eqref{eq:linear-Gaussian-model}.

\begin{remark}\label{rem:KF_derivation}
The solution of the optimal control problem yields the optimal control
input $u^\opt=\{u^\opt_t:0\leq t\leq T\}$, along with the vector $y_0\in\Re^d$ that determines the
minimum-variance estimator:
\begin{align*}
	S_T &= \mu(y_0^\tp x) - \int_0^T \big(u_t^\opt \big)^\tp \ud Z_t
	= y_0^\tp m_0 - \int_0^T \big(u_t^\opt \big)^\tp \ud Z_t
\end{align*}
The Kalman filter is obtained by expressing $\{S_t(f) : t\ge 0,\
f\in\Re^d\}$ as the solution to a linear SDE~\cite[Ch.~7.6]{astrom1970}.
\end{remark}

\section{Solution of the optimal control problem}\label{sec:standard-form}

The BSDE constrained optimal control problem~\eqref{eq:dual-optimal-control} is not in its standard form~\cite[Eq. 5.10]{pardoux2014stochastic}. There are two issues:
\begin{itemize}
	\item {\bf The probability space:} The driving martingale of the BSDE~\eqref{eq:dual-optimal-control-b} is $Z$, which is a $\tsP$-B.M. However, the expectation in defining the optimal control objective~\eqref{eq:dual-optimal-control-a} is with respect to the measure $\sP$. 
	\item {\bf The filtration:} The `state' of the optimal control problem $(Y,V)$ is adapted to the filtration $\clZ$. However, the cost function~\eqref{eq:dual-optimal-control-a} also depends upon the non-adapted exogenous process $X$.  
\end{itemize}

The second problem is easily fixed by using the tower property of conditional expectation.
To resolve the first problem, we have two choices:
\begin{enumerate}
	\item Use the change of measure to evaluate $\sJ_T(U)$ with respect to $\tsP$ measure, or
	\item Express the BSDE using a driving martingale that is a $\sP$-B.M. A convenient such process is the innovation process. 
\end{enumerate}
In this paper, the standard form of the dual optimal control problem
is presented based on the first choice.  For an analysis based on the
second choice, see~\cite{kim2019duality} and~\cite[Sec.~5.5]{JinPhDthesis}.    

\medskip

In order to express the expectation for the control
objective~\eqref{eq:dual-optimal-control-a} with respect to $\tsP$, we
use the change of measure (see Appendix~\ref{ss:derivation-lagrangian} for the calculation) to obtain
\begin{align*}
	\sJ_T(U) &= \var_0(Y_0) + \tE\Big( \int_0^T \ell(Y_t,V_t,U_t;\sigma_t)\ud t\Big)
\end{align*}
where the \emph{Lagrangian} $\ell:\clY\times \clY^m\times \Re^m \times \clY^\dagger\to \Re$ is defined by
\begin{equation}
	\ell(y,v,u;\rho) = \rho\big(\Gamma y \big) + \rho\big(|u+v|^2\big)
\end{equation}

The dual optimal control problem (standard form) is now expressed as follows:	
\begin{subequations}\label{eq:dual-optimal-control-std}
	\begin{align}
&		\mathop{\text{Minimize}}_{U\in\clU}\quad \sJ_T(U)
          = \dvar_0(Y_0) + \tE\Big(  \int_0^T \ell(Y_t,V_t,U_t;\sigma_t)\ud t
		\Big) \label{eq:dual-optimal-control-std-a}\\
&		\text{\rm Subject to:}\nonumber \\ 
&  -\! \ud Y_t(x) = \big((\clA Y_t)(x) + h^\tp(x)(U_t+V_t(x))\big)\ud t - V_t^\tp(x)\ud Z_t\nonumber \\
&    \quad \;  Y_T (x)  = F(x), \;\; x \in \bS 
\label{eq:dual-optimal-control-std-b}
	\end{align}
\end{subequations}


\medskip

\begin{remark}
	The Lagrangian is a time-dependent random functional of the dual state $(y,v)$ and the control $u$. The randomness and time-dependency comes only from the last argument $\sigma_t$.
\end{remark}

\medskip

\subsection{Solution using the maximum principle}

Because $y\in\clY$ is a function, the co-state $p\in \clY^\dagger$ is
a measure.  The \emph{Hamiltonian} $\clH: \clY \times \clY^m\times 
\Re^m\times \clY^\dagger\times \clY^\dagger \to \Re$ is defined as follows:
\[
\clH(y,v,u,p;\rho) = -p\big(\clA y + h^\tp (u + v) \big)- \ell(y,v,u;\rho)
\]
In the following, the Hamilton's equations for the optimal trajectory are derived by 
an application of the maximum principle for BSDE constrained optimal
control problems~\cite[Thm.~4.4]{peng1993backward}. 

The Hamilton's equations are expressed in terms of the 
derivatives of the Hamiltonian.  In order to take derivatives with
respect to functions and measures, we  adopt the notion of G\^ateaux
differentiability. 
Given a nonlinear functional
$F:\clY\to \Re$, 
the G\^ateaux derivative $F_y(y)\in \clY^\dagger$ is obtained from the
defining relation~\cite[Sec.~10.1.3]{bensoussan2018estimation}:
\[
\frac{\ud}{\ud \varepsilon} F(y + \varepsilon\tilde{y})\Big|_{\varepsilon=0} = 
\big\langle \tilde{y}, F_y(y)\big\rangle,\quad \forall\,\tilde{y} \in \clY
\]

For the problem at hand, the derivatives of the 
Hamiltonian are as follows:
\begin{align*}
	\clH_y(y,v,u,p;\rho) &= -\clA^\dagger p - \big(\rho(\Gamma y)\big)_y\\
	\clH_v(y,v,u,p;\rho) &= -ph - 2(u+v)\rho\\
	\clH_u(y,v,u,p;\rho) &= -p(h) - 2\rho(\ones)u - 2\rho(v)\\
	\clH_p(y,v,u,p;\rho) &= -\clA y - h^\tp(u+v)
\end{align*}
where $\clA^\dagger$ is the adjoint of $\clA$ (whereby
$(\clA^\dagger\rho) (f) = \rho(\clA f)$ for all $f\in\clY,\rho\in\clY^\dagger$).  
Using this notation, the Hamilton's equations are as follows: 

\medskip

\begin{theorem}\label{thm:optimal-solution}
	Consider the optimal control problem~\eqref{eq:dual-optimal-control-std}.  Suppose
	$U^\opt$ is the optimal control input and the $(Y,V)$ is the
        associated solution of the BSDE~\eqref{eq:dual-optimal-control-std-b}.
	Then there exists a $\clZ$-adapted measure-valued stochastic process $P=\{P_t:0\le t \le T\}$ such
	that
	\begin{subequations}\label{eq:Hamilton_eqns}
		\begin{flalign}
&			 \ud P_t = -\clH_y(Y_t,
V_t,U_t^\opt,P_t;\sigma_t)\ud t  \nonumber \\ &\qquad \qquad - \clH_v^\tp(Y_t, V_t,U_t^\opt,P_t;\sigma_t) \ud Z_t
			\label{eq:Hamilton_eqns2-a}\\
&			\ud Y_t = \clH_p(Y_t,V_t,U_t^\opt,P_t;\sigma_t) \ud t + V_t\ud Z_t \label{eq:Hamilton_eqns2-b}\\
&			 \frac{\ud P_0}{\ud \mu}(x)= 2\big(Y_0(x)-\mu(Y_0)\big),\quad Y_T(x) = F(x),\quad  x\in\bS  \label{eq:Hamilton_eqns2-c}
		\end{flalign}
	\end{subequations}
	where the optimal control is given by
	\begin{equation}\label{eq:opt-cont-soln}
		U_t^\opt = -\frac{1}{2} \frac{P_t(h)}{\sigma_t(\ones)} - \pi_t(V_t),\quad \tsP\text{-a.s.},\; 0\le t \le T
	\end{equation}
(In~\eqref{eq:Hamilton_eqns2-c}, $\frac{\ud P_0}{\ud \mu}$ denotes the R-N derivative of the
measure $P_0$ with respect to the measure $\mu_0$). 
\end{theorem}

\medskip

\begin{proof}
See Appendix~\ref{ss:pf-thm52}.
\end{proof}

\medskip

\begin{remark}\label{rem:state-costate-linear-trans}
	From linear optimal control theory, it is known that $P_t$ is
        related to $Y_t$ by a ($\clZ_t$-measurable) linear
        transformation~\cite[Sec.~6.6]{yong1999stochastic}. The
        boundary condition $\frac{\ud P_0}{\ud \mu}(x) =
        2\big(Y_0(x)-\mu(Y_0)\big)$ suggests that the R-N derivative 
	\begin{equation}\label{eq:P-t-ansatz}
		\frac{\ud P_t}{\ud \sigma_t}(x) = 2\big(Y_t(x) - \pi_t(Y_t)\big),\quad \tsP\text{-a.s.},\; 0\le t \le T,\; x\in \bS 
	\end{equation}
	This is indeed the case as we show in
        Appendix~\ref{ssec:pf-ansatz} by verifying that~\eqref{eq:P-t-ansatz} solves the Hamilton's equation. 
	Combining this formula with~\eqref{eq:opt-cont-soln}, we have a formula for optimal control input as a feedback control law:
	\begin{equation*}
		U_t^\opt = -\big(\pi_t(hY_t) - \pi_t(h)\pi_t(Y_t)\big) - \pi_t(V_t),\quad 0\le t \le T
	\end{equation*}
\end{remark}

\medskip
\medskip

\subsection{Explicit formulae for the guiding examples}

\begin{example}[Finite state-space] (Continued from Example~\ref{ex:finite}).  A real-valued function $f$
  (resp. a measure $\rho$) is identified with a column vector in
  $\Re^d$ where the $i^{\text{th}}$ element of the vector represents
  $f(i)$ (resp. $\rho(i)$), and $\rho(f) = \rho^\tp f$.  In this
  manner, the generator $\clA$ is identified with a rate matrix
  $A\in\Re^{d\times d}$ and the 
observation function $h$ is identified with a matrix
$H\in\Re^{d\times m}$. Let $\{e_1,e_2,\hdots,e_d\}$ denote the
canonical basis in $\Re^d$, $Q(i) = \sum_{j \in \bS}
A(i,j)(e_i-e_j)(e_i-e_j)^\tp$ and $\rho(Q) = \sum_{i \in
  \bS}\rho(i)Q(i)$.  For any vector $b\in\Re^d$, $B=\dv(b)$ is a
$d\times d$ diagonal matrix whose diagonal entires are defined as
$B(i,i) = b(i)$ for $i=1,2,\hdots,d$.  For a $d\times d$ matrix $B$,  $b=\dv^\dagger(B)$ is
a $d$-dimensional vector whose entries are defined as $b(i) =
B(i,i)$ for $i=1,2,\hdots,d$.

The Lagrangian $\ell:\Re^d\times \Re^{d\times m} \times \Re^m \times
\Re^d \to \Re$ and the Hamiltonian $\clH:\Re^d\times \Re^{d\times m}
\times \Re^m \times \Re^d\times \Re^d \to \Re$ are as follows:
\begin{align*}
&\ell(y,v,u;\rho) = y^\tp \rho(Q)y + \rho(\ones)|u|^2 + 2u^\tp v\rho +
  \rho^\tp\dv^\dagger (vv^\tp)\\
&\clH(y,v,u,p;\rho) = -p^\tp (Ay +Hu + \dv^\dagger(Hv^\tp))- \ell(y,v,u;\rho) 
\end{align*}
The functional derivatives are now the partial derivatives.  For the
Hamiltonian, these are as follows:
\begin{align*}
	\clH_y(y,v,u,p;\rho) &= -A^\tp p - 2\rho(Q)y\\
	\clH_v(y,v,u,p;\rho) &= -\dv(p)H - 2\rho u^\tp - 2\dv(\rho) v\\
	\clH_u(y,v,u,p;\rho) &= -H^\tp p - 2\rho(\ones)u -2v^\tp \rho\\
	\clH_p(y,v,u,p;\rho) &= -A y - Hu - \dv^\dagger(Hv^\tp)
\end{align*}
The Hamilton's equations are given by	
\begin{flalign*}
	\ud P_t &= \big(A^\tp P_t + 2\sigma_t(Q)Y_t\big) \ud t
        \\\nonumber
&\qquad + \big(\dv(P_t)H + 2\sigma_tU_t^\tp +2\dv(\sigma_t)V_t\big) \ud Z_t
	\\
	 \ud Y_t &= -\big(AY_t+HU_t+\dv^\dagger(HV_t^\tp)\big) \ud t + V_t \ud Z_t \\
	P_0 & = 2\Sigma_0 Y_0,\quad Y_T = F 
\end{flalign*}
where $\Sigma_0 := \dv(\mu) - \mu \mu^\top$. 
\end{example}

\medskip

\begin{example}[Euclidean state-space]  (Continued from Example~\ref{ex:itodiffusion}).  We consider the It\^{o}
  diffusion~\eqref{eq:dyn_sde} in $\Re^d$ with a prior density
  denoted as $\mu$.  Likewise, $\rho$ and $p$ are 
  used to denote the density of the respective measures.  Doing so,
  the Lagrangian and the Hamiltonian are as follows:  
\begin{align*}
\ell(y,v,u;\rho) & = \int_{\Re^d} \rho(x)\big(|\sigma^\tp(x) \nabla
  y(x)|^2 + |u+v(x)|^2\big)\ud x\\
\clH(y,v,u,p;\rho) &= -\int_{\Re^d} p(x)\big(\clA y(x) +
  h^\tp(x)(u+v(x))\big)\ud x \\
&\qquad \qquad-\ell(y,v,u;\rho)
\end{align*}
The functional derivatives are computed by evaluating the first
variation. These are as follows:
\begin{align*}
	\clH_y(y,v,u,p;\rho) &= -\clA^\dagger p + 2\divg\big(\sigma\sigma^\tp(\nabla y) \rho\big)\\
	\clH_v(y,v,u,p;\rho) &= -ph-2(u+v)\rho\\
	\clH_u(y,v,u,p;\rho) &= -p(h) - 2\rho(\ones)u - 2\rho(v)\\
	\clH_p(y,v,u,p;\rho) &= -\clA y - h^\tp(u+v)
\end{align*}
where $\rho(v)$ is now understood to mean $\int\rho(x)v(x) \ud x$ and the
formula for adjoint is
\[
({\cal A}^\dagger p)(x) = -\divg(a p)(x) + \half  \sum_{i,j=1}^d \frac{\partial^2}{\partial x_i \partial x_j}\big([\sigma\sigma^\tp]_{ij} p\big)(x)
\]
Therefore, the Hamilton's equations are given by	
\begin{flalign*}
\ud P_t(x) &= \big((\clA^\dagger P_t)(x) - 2\divg \big(\sigma\sigma^\tp
(\nabla Y_t) \sigma_t\big)(x)\big) \ud t 
\\ & \qquad + \big(P_t(x)h(x) + 2(U_t+V_t(x))\sigma_t(x)\big) \ud Z_t
	 \\
\ud Y_t(x) &= -\big(\clA Y_t+h^\tp(x)(U_t+V_t(x))\big) \ud t +
V_t^\tp(x) \ud Z_t 
\\
P_0(x) &= 2\mu(x)\big(Y_0(x)-\mu(Y_0)\big),\; Y_T(x) = F(x),\; x\in \Re^d
\end{flalign*}
where note that $P_t$ is now a (random) function (same as $Y_t$). 
\end{example}

\medskip

\section{Martingale characterization}\label{ssec:martingale}

Although we do not have an HJB equation, a martingale characterization
of the optimal solution is possible as described in the following theorem:


\begin{theorem}\label{thm:martingale}
	Fix $U\in \clU$. Consider a $\clZ$-adapted real-valued
        stochastic process $M = \{M_t:0\le t \le T\}$
	\begin{equation*}\label{eq:martingale-characterization}
		M_t := \clV_t(Y_t) - \int_0^t \ell(Y_s,V_s,U_s;\pi_s)\ud s,\quad 0\le t\le T
	\end{equation*}
	where $(Y,V)$ is the solution to the
        BSDE~\eqref{eq:dual-optimal-control-b} and $\pi$ is the nonlinear filter.
	Then $M$ is a $\sP$-supermartingale, and $M$ is a $\sP$-martingale if and only if
	\begin{equation}\label{eq:optimal-solution}
		U_t = -\big(\pi_t(hY_t) - \pi_t(h)\pi_t(Y_t)\big) - \pi_t(V_t)
	\end{equation}
for $ 0\le t \le T$, $\sP\text{-a.s.}$.
\end{theorem} 

\medskip

\begin{proof}
See Appendix~\ref{ss:pf-thm53}.
\end{proof}

\medskip

A direct consequence of~\Thm{thm:martingale} is the optimality of the control~\eqref{eq:optimal-solution}, because
\[
\E(M_T) \le \E(M_0)
\]
which means
\[
\E\big(\clV_T(F)\big) \le \E\Big(\clV_0(Y_0) +
\int_0^T\ell(Y_t,V_t,U_t;\pi_t)\ud t\Big) = \sJ_T(U)
\]
with equality if and only if $U = U^\opt$.  Therefore, the expected
value of the conditional variance $\var_T(F)=\E\big(\clV_T(F)\big)$ is
the optimal value functional for the optimal control problem.  
%
%

\medskip

\begin{remark}
We now have a complete solution of the optimal control
problem~\eqref{eq:dual-optimal-control}.  Remarkably, the solution
admits a meaningful interpretation not only at the terminal time $T$
but also for intermediate times $0\leq t\leq T$.     
At time $t$,
\begin{itemize}
\item The optimal value functional is $\var_t(Y_t)$ (formula~\eqref{eq:estimator-t-variance}).
\item The optimal control $U_t^\opt$ is a feedback control law~\eqref{eq:optimal-solution}.  
\item The optimal estimate is $\pi_t(Y_t)$ (formula~\eqref{eq:estimator-t}).
\end{itemize}
Formula~\eqref{eq:estimator-t} for $\pi_t(Y_t)$ explicitly connects the
optimal control to the optimal filter. In particular, the optimal
control up to time $t$ yields an optimal estimate of $Y_t(X_t)$.    


Because of the BSDE constrained nature of the optimal control
problem~\eqref{eq:dual-optimal-control}, an explicit characterization of
the optimal value functional and the feedback
form of the optimal control are both
welcome surprises.  It is noted that the feedback formula~\eqref{eq:optimal-solution} for the
optimal control is derived using two approaches: using the maximum principle
(Rem.~\ref{rem:state-costate-linear-trans}) and using the martingale
characterization (\Thm{thm:martingale}). 
\end{remark}

\section{Derivation of the nonlinear filter}\label{sec:derivation-of-the-filter}

From~\Prop{prop:optimal-solution}, using the formula~\eqref{eq:optimal-solution} for optimal
control
\begin{equation}\label{eq:estimator-with-optimal-ctrl}
	\pi_t(Y_t) = \mu(Y_0) + \int_0^t \big(\pi_t(hY_s) - \pi_s(h)\pi_s(Y_s)+\pi_s(V_s)\big)^\tp\ud Z_s
\end{equation}
for $ 0\le t \le T$, $\sP\text{-a.s.}$.  Because the equation for 
$Y$ is known, a natural question is whether~\eqref{eq:estimator-with-optimal-ctrl} can be used
to obtain the equation for nonlinear filter (akin to the derivation of
the Kalman filter described in Rem.~\ref{rem:KF_derivation}).  A formal
derivation of the nonlinear filter along these lines is given in
Appendix~\ref{ssec:pf-nonlinear-filter-derivation}.


	



\begin{table*}
	\centering
	\renewcommand{\arraystretch}{2}
	\caption{Comparison of the Mitter-Newton duality and the
          duality proposed in this paper} \label{tb:comparison}
	\small
	\begin{tabular}{p{0.25\textwidth}p{0.35\textwidth}p{0.35\textwidth}}
		& {\bf Mitter-Newton duality} & {\bf Duality proposed
                                                in this paper} \\
          \hline \hline
		Filtering/smoothing objective & Minimize relative
                                                entropy (Eq.~\eqref{eq:rel-entropy-mitter})
                                              & Minimize
                                                variance (Eq.~\eqref{eq:minimum-variance})
          \\ \hline
		Observation (output) process & Pathwise ($z$ is a sample path) & $Z$ is a
          stochastic process\\ \hline
		Control (input) process & $U_t$ has the dimension of
                                        the process noise & $U$ and
                                                            $Z$ are
                                                            both
                                                            elements
                                                            of $L^2_{\clZ}([0,T];\Re^m)$
                                                        \\ \hline
		Dual optimal control problem &
                                               Eq.~\eqref{eq:opt-cont-sde-hjb-intro} 
                                              &  Eq.~\eqref{eq:dual-optimal-control}
          \\ \hline
		Arrow of time & Forward in time & Backward in time \\
          \hline  
		Dual state-space & $\bS$: same as the state-space for $X_t$ &
                                                                      $\clY$:
                                                                             the
                                                                             space
                                                                             of functions on $\bS$ \\ \hline 
		Constraint & Controlled copy of the state
                             process SDE~\eqref{eq:opt-cont-sde-hjb-intro-a}
                                              & Dual control system
                                                BSDE~\eqref{eq:dual-optimal-control-b}
          \\ \hline
Running cost (Lagrangian) & $l(x,u\,;z_t)= \half |u|^2 + \half h^2(x) + z_t(\clA^u h)(x)$ &
$l(y,v,u;x) = (\Gamma y)(x) + |u+v(x)|^2$ \\ \hline 
Value function (its interpretation) & Minus log of the posterior density & Expected value of the conditional variance  \\
          \hline  
Asymptotic analysis (condition) & Unclear & Stabilizability of BSDE $\Leftrightarrow$  Detectability of HMM\\
          \hline  
Optimal solution gives & Forward-backward equations of smoothing  & 
 Equation of nonlinear filtering \\
          \hline 
		Linear-Gaussian special case & Minimum energy duality~\eqref{eq:mee-intro} & Minimum variance duality~\eqref{eq:mv-intro}  \\ \hline 
\hline
	\end{tabular}
\end{table*}
\section{Comparison with Mitter-Newton Duality}\label{sec:comp}


\subsection{Review of Mitter-Newton duality}

In~\cite{mitter2003}, Mitter and Newton introduced a modified version
of the Markov process $X$.  The modified process is denoted by $\tilde{X} :=\{\tilde{X}_t:0\le t\le
T\}$.  The problem is to pick (i) the initial prior
$\tilde{\mu}$; and (ii) the state transition, such that the
probability law of
$\tilde{X}$ equals the conditional law for $X$.  

This is accomplished by setting up an optimization problem on the
space of probability laws. 
Let $\sP_X$ denote the law for $X$, $\sQ$ denote the law for
$\tilde{X}$, and $\sP_{X\mid z}$ denote the conditional law for $X$
given an observation sample path
$z=\{z_t\in\Re^m:0\le t\le T\}$.  Assuming
$\sQ\ll\sP_X$, the objective function is the relative entropy between $\sQ $ and $\sP_{X\mid z}$:
\begin{equation}\label{eq:rel-entropy-mitter}
	\min_{\sQ} \quad \E_{\sQ}\Big(\log \frac{\ud \sQ}{\ud \sP_X}\Big) - \E_{\sQ}\Big(\log\frac{\ud \sP_{X\mid z}}{\ud \sP_X}\Big) 
\end{equation}
In~\cite{van2006filtering},~\eqref{eq:rel-entropy-mitter} is referred to as the variational Kallianpur-Striebel formula.  
For Example~\ref{ex:itodiffusion} (It\^o diffusion), this procedure yields the following
stochastic optimal control problem:
\begin{subequations}\label{eq:opt-cont-sde-hjb-intro}
	\begin{align}
	&\mathop{\text{Min}}_{\tilde{\mu}, \; U}: \;\; \sJ(\tilde{\mu},U\,;z)  \nonumber\\
	&\qquad= \E\Big(\log \frac{\ud \tilde{\mu}}{\ud \mu}(\tilde{X}_0) - z_T h(\tilde{X}_T)
  + \int_0^T l(\tilde{X}_t,U_t\,;z_t)\ud t\Big) \label{eq:opt-cont-sde-hjb-intro-a}\\
&		\text{Subj.} : \;\; \ud \tilde{X}_t = a(\tilde{X}_t)\ud t +
		\sigma(\tilde{X}_t)(U_t\ud t +
		\ud \tilde{B}_t), \;\;
		\tilde{X}_0 \sim \tilde{\mu} \label{eq:opt-cont-sde-hjb-intro-b}
	\end{align}
\end{subequations}
where 
\begin{equation*}\label{eq:lagrangian-mitter}
l(x,u\,;z_t) := \half |u|^2 + \half h^2(x) + z_t(\clA^u h)(x)
\end{equation*}
where
$\clA^u$ is the generator of the controlled Markov process
$\tilde{X}$.  A similar construction is also possible for 
Example~\ref{ex:finite} (finite
state-space)~\cite[Sec.~2.2.2]{van2006filtering},~\cite[Sec.~3.3]{kim2020smoothing}.




The problem~\eqref{eq:opt-cont-sde-hjb-intro} is a standard stochastic
optimal control problem whose solution is obtained by writing the 
HJB equation (see~\cite{kim2020smoothing}),
\begin{align*}
-\frac{\partial v_t}{\partial t}(x) &= \big({\cal A}(v_t+z_th)\big)(x)
                                      + \half h^2(x) \\
& \qquad -\half|\sigma^\tp\nabla (v_t+z_th)(x)|^2\\
v_T(x) &= - z_Th(x),\quad x\in\Re^d
\end{align*}
and the optimal control $
	U_t = u_t^\opt(\tilde{X}_t) 
$
	where 
\[
u_t^\opt(x) = -\sigma^\tp \nabla(v_t + z_th)(x)
\]
By expressing the value function 
\[
v_t(x) = -\log \big(q_t(x)e^{z_th(x)}\big)
\]
a direct calculation shows that the process $\{q_t:0\leq t\leq T\}$ satisfies the
backward Zakai equation of the smoothing
problem~\cite{pardoux1979backward},\cite[Thm.~3.8]{pardoux1981non}.  This shows the connection to both the log transformation and
to the smoothing problem.  In fact, the above can be used to derive
the forward-backward equations of nonlinear
smoothing~(see \cite{kim2020smoothing} and~\cite[Appdx.~B]{JinPhDthesis}).

\medskip

\begin{remark}\label{rm:density-dynamics}
The stochastic optimal control
problem~\eqref{eq:opt-cont-sde-hjb-intro} is equivalently stated as a
deterministic optimal control problem on
$\clY^\dagger$~\cite[Sec.~3.2]{kim2020smoothing}. 
Note that the optimal control problem depends on a (fixed)
observation sample path $z$, which is the reason why a deterministic
formulation is available.
\end{remark}

%
%

\subsection{Linear Gaussian case}\label{sec:MN-LG}

The goal is to relate~\eqref{eq:opt-cont-sde-hjb-intro} to the minimum
energy duality~\eqref{eq:mee-intro} described
in~\Sec{sec:introduction} for the linear Gaussian
model~\eqref{eq:linear-Gaussian-model}.  In the linear Gaussian case, the controlled
process~\eqref{eq:opt-cont-sde-hjb-intro-b} becomes
\begin{equation}\label{eq:tildeX-MN-G}
\ud \tilde{X}_t = A^\tp \tilde{X}_t \ud t + \sigma U_t \ud t + \sigma \ud \tilde{B}_t,\quad \tilde{X}_0\sim N(\tilde{m}_0,\tilde{\Sigma}_0)
\end{equation}
where $U$, $\tilde{m}_0,\tilde{\Sigma}_0$ are decision variables.
Because the problem is linear Gaussian, it suffices to consider a
linear control law of the form
\begin{equation}\label{eq:U-MN-G}
U_t = K_t(\tilde{X}_t-\tilde{m}_t) + u_t
\end{equation}
where $\tilde{m}_t := \E(\tilde{X}_t)$ and the two deterministic
processes
\begin{align*}
K&=\{K_t\in \Re^{p\times
  d}:0\leq t\leq T\}\\
u&=\{u_t\in \Re^p:0\leq t\leq T\}
\end{align*}
are the new decision variables.  With a linear control law~\eqref{eq:U-MN-G}, the state
$\tilde{X}_t$ is a Gaussian random variable with mean $\tilde{m}_t$
and variance $\tilde{\Sigma}_t$.  It is possible to equivalently express~\eqref{eq:opt-cont-sde-hjb-intro}
as two un-coupled deterministic optimal control problems, for the mean
and for the variance, respectively.  Detailed calculations showing
this are contained in Appendix~\ref{apdx:Mitter-LG-case}.  In particular, it is
shown that the optimal control problem for the mean is the
classical minimum
energy duality~\eqref{eq:mee-intro}.  

\subsection{Comparison}

Table~\ref{tb:comparison} provides a side-by-side comparison of the
two types of duality:
\begin{itemize}
\item Mitter-Newton duality~\eqref{eq:opt-cont-sde-hjb-intro} on the
  left-hand side; and
\item Duality~\eqref{eq:dual-optimal-control} proposed in this paper on the right-hand side.
\end{itemize}
In \Sec{sec:MN-LG} and~\Sec{sec:LG-case}, the two are shown to be
generalization of the classical minimum energy
duality~\eqref{eq:mee-intro} and the minimum variance
duality~\eqref{eq:mv-intro}, respectively.  All of this conclusively
answers the two questions
raised in~\Sec{sec:introduction}.  

We make a note of some important distinctions (compare with the bulleted
list in~\Sec{sec:introduction}):
\begin{itemize}
\item{\em Inputs and outputs.} In proposed duality~\eqref{eq:dual-optimal-control}, inputs and outputs
  are dual processes that have the same dimension.  These are element of
  the same Hilbert space $\clU$.    
\item{\em Constraint.} The constraint is the dual control system~\eqref{eq:dual-optimal-control-b}
  studied in the companion paper (part I).  
\item{\em Stability condition.} For asymptotic analysis
  of~\eqref{eq:dual-optimal-control}, stabilizability of the
  constraint is the most natural condition.  The main result of part I was
  to establish that 
  stabilizability of the dual control system is equivalent to the
  detectability of the HMM.  The latter condition of course is central to filter
  stability.  
\item{\em Arrow of time.} The dual control system is backward in
  time.  However, it is important to note that the information
  structure (filtration) is forward in time.  In particular, all the
  processes are forward adapted to the filtration $\clZ$ defined by the
  observation process.  
\end{itemize}   
A major drawback of the proposed duality is that the problem (for
the Euclidean state-space $\bS=\Re^d$) is infinite-dimensional.  This is to be expected because the nonlinear filter is
infinite-dimensional.
In contrast, the state space in the minimum energy duality is $\Re^d$ which is important for algorithm design as in MEE.
Having said that, the linear quadratic nature of the infinite-dimensional problem may
prove to be useful in practical applications of this work. 

\section{Conclusions and directions of future work}
\label{sec:conc}

In this paper, we presented the minimum variance dual optimal control 
problem for the nonlinear filtering problem.  The mathematical
relationship between the two problems is given by a duality
principle.  Two approaches are described to solve the problem, based
on maximum principle and based on a martingale characterization.  A
formula for the optimal control as a feedback control law is obtained,
and used to derive the equation of the nonlinear filter.  A detailed
comparison with the Mitter-Newton duality is given.  

There are several possible directions of future research:  An 
important next step is to use the controllability and stabilizability
definitions of the dual control system to recover the known results in
filter stability.  Research on this has already begun with preliminary
results appearing in~\cite[Chapter 7-8]{JinPhDthesis}
and~\cite{kim2021ergodic,kim2021detectable}.  Although some sufficient
conditions have been obtained and compared with literature, a
complete resolution still remains open.  

Both the stability analysis and the optimal control formulation
suggest natural connections to the dissipativity theory.  Because
the dual control system is linear, one might consider quadratic forms
of supply rate function as follows (compare with the formula for the
running cost $l$):
\[
s(y, v, u; x) := \gamma |u + v(x)|^2 - |y(x) - {c}_t|^2
\] 
where $\gamma>0$ and $c:=\{c_t:0\leq t\leq T\} \in
L^2_{\clZ}\big([0,T];\Re \big)$ is a suitable stochastic process
(which can be picked).
Establishing conditions for existence of a storage function and
relating these conditions to the properties of the HMM may be useful
for stability and robustness analysis.


Another avenue is numerical approximation of the nonlinear filter by
considering sub-optimal solutions of the dual optimal control problem.
The simplest choice is to consider deterministic control inputs $U \in
L^2\big([0,T];\Re^m\big)$. 
Some preliminary work on algorithm design along these lines appears
in~\cite[Rem.~1]{kim2019duality},~\cite[Sec.~9.2]{JinPhDthesis}
and~\cite[Ch.~4]{jan2021master}.  
In particular for the finite state space case, this approach provides derivation and justification of Kalman filter for Markov chains~\cite{lipkrirub84}.
In this regard, it is useful to relate duality to both the feedback particle filter
(FPF)~\cite{taoyang_TAC12} and to the special cases (apart from the
linear Gaussian case) where the optimal filter is known to be
finite-dimensional, e.g.~\cite{benevs1981exact}.



\section{Acknowledgement}

It is a pleasure to acknowledge Sean Meyn and Amirhossein Taghvaei for many useful
technical discussions over the years on the topic of duality.  The
authors also acknowledge Alain Bensoussan for his early encouragement of this work.  

\appendix

\section{Proofs of the statements}

\subsection{Proof of~\Thm{thm:duality-principle}}\label{ss:pf-thm51}

For a Markov process, the following process is a martingale:
\begin{equation*}
	N_t(g) = g(X_t) - \int_0^t \clA g(X_s)\ud s
\end{equation*}
Upon applying the It\^o-Wentzell
theorem~\cite[Thm.~1.17]{rozovskiui2018stochastic} on $Y_t(X_t)$ (note here
that all stochastic processes
are forward adapted),
\begin{align*}
	\ud Y_t(X_t) 
	&= - U_t^\tp \ud Z_t + \big(U_t + V_t(X_t)\big)^\tp \ud W_t + \ud N_t(Y_t)
\end{align*}
Integrating both sides from $0$ to $T$,
\begin{align*}
F(X_T) &= Y_0(X_0) - \int_0^T U_t^\tp \ud Z_t \\
&+ \int_0^T (U_t+V_t(X_t))^\tp \ud W_t + \int_0^T \ud N_t(Y_t)
\end{align*}
Consider now an estimator
\begin{equation*}\label{eq:estimator-general}
	S_T = b- \int_0^T U_t^\tp \ud Z_t
\end{equation*}
where $b\in \Re$ is a deterministic constant. Then
\begin{align*}
F(X_T) - S_T &= \big(Y_0(X_0)-b\big) + \int_0^T (U_t+V_t(X_t))^\tp \ud W_t \\
&\quad +\int_0^T \ud N_t(Y_t)
\end{align*}
The left-hand side is the error of the
estimator. The three terms on the
right-hand side are mutually independent. Therefore, upon squaring and
taking an expectation
\begin{align*}
\E\big(|F(X_T)-S_T|^2\big)= \E\big(|Y_0(X_0)-\mu(Y_0)|^2\big) + (\mu(Y_0)-b)^2\\
+ \E\Big(\int_0^T |U_t+V_t(X_t)|^2 + (\Gamma Y_t)(X_t) \ud t \Big)
\end{align*}
The proof is completed by setting $b = \mu(Y_0)$.

\subsection{Proof of~\Lemma{lem:representation}}\label{ss:pf-lem51}

Because $Z$ is a $\tsP$-B.M., the formula holds for $\pi_T(F)\in
L^2_{\clZ_T}(\Omega;\Re)$ by the Brownian motion representation
theorem~\cite[Thm.~5.18]{le2016brownian}. Note that 
\[
|\pi_T(F)|^2 \le \|F\|_\clY^2,\quad \tsP\text{-a.s.}
\]
because $\|\cdot\|_\clY$ is the sup norm.  
Therefore if $F\in L^2_{\clZ_T}(\Omega;\clY)$ then $\pi_T(F)\in
L^2_{\clZ_T}(\Omega;\Re)$. The conclusion follows.


\subsection{Proof of~\Prop{prop:optimal-solution}}\label{ssec:pf-optimal-solution}

Using optimal control 
$U^\opt =\{U^\opt_t:0\leq t\leq T\}\in \clU$,
$(Y,V) =\{(Y_t,V_t):0\leq t\leq T\} \in
L^2_{\clZ}\big([0,T];\clY\times \clY^m\big)$ is the solution of the
BSDE~\eqref{eq:dual-optimal-control-b} with $Y_T = F \in L^2_{\clZ_T}(\Omega;\clY)$.   Fix $t \in [0,T]$ and let
\[
S_t = \mu(Y_0) - \int_0^t \big(U_s^\opt\big)^\tp \ud Z_s
\]
Then by repeating the proof of~\Thm{thm:duality-principle} now over the
time-horizon $[0,t]$,
\[
\E\big(|Y_t(X_t)-S_t|^2\big) = \var_0(Y_0) + \E (\int_0^t l(Y_s,V_s,U_s^\opt;X_s) \ud s)
\]
If $\E\big(|Y_t(X_t)-S_t|^2\big) =\var_t(Y_t)$ then there is nothing
to prove.  Because then $S_t = \pi_t(Y_t)$ ($\sP$-a.s.) by the
uniqueness of the conditional expectation.  Therefore, suppose 
\[
\var_t(Y_t)=\E \big(|Y_t(X_t)-\pi_t(Y_t)|^2\big) < \E \big(|Y_t(X_t)-S_t|^2\big) 
\]
In this case, we show that there exists a $\tilde{U} \in \clU$ such
that $\sJ_T(\tilde{U}) < \sJ_T(U^\opt)$.  Because $U^\opt$ is the
optimal control, this provides the necessary contradiction. 

Set $C:= \E(\int_t^T 
l(Y_s,V_s,U_s^\opt;X_s) \ud s)$ and we have
\[
\sJ_T(U^\opt) = \E\big(|Y_t(X_t)-S_t|^2\big)  + C
\]
Because $Y_t \in L^2_{\clZ_t}(\Omega;\clY)$,
by~\Lemma{lem:representation} there exists $\hat{U}\in
L_{\clZ}^2([0,t];\Re^m)$ such that
\[
\pi_t(Y_t) = \tE(\pi_t(Y_t)) - \int_0^t \hat{U}_s^\tp \ud Z_s,\quad \tsP\text{-a.s.}
\]
Consider an admissible control $\tilde{U}$ as follows
\[
\tilde{U}_s = \begin{cases}
	\hat{U}_s & s \le t\\
	U_s^\opt &  s > t
\end{cases}
\]
and denote by $(\tilde{Y},\tilde{V})$ the solution of the BSDE with
the control $\tilde{U}$.  Because of the uniqueness of the solution,
$(\tilde{Y}_s,\tilde{V}_s) = (Y_s,V_s)$ for all $s>  t$ and therefore
\begin{align*}
	\sJ_T(\tilde{U}) &= \E \big(|Y_t(X_t)-\pi_t(Y_t)|^2\big) + C\\
	&< \E \big(|Y_t(X_t)-S_t|^2\big)  + C = \sJ_T(U^\opt)
\end{align*}
This supplies the necessary contradiction and completes the proof. 

\subsection{Derivation of the Lagrangian}\label{ss:derivation-lagrangian}

Using the change of measure formula~\eqref{eq:COM_formula},
\begin{align*}
\E \big((\Gamma Y_t)(X_t)\big) &= \tE \big(\sigma_t(\Gamma Y_t)\big) \\
\E\big(|U_t +
V_t(X_t)|^2\big) & = \tE\big( \sigma_t(|U_t + V_t|^2)\big)
\end{align*}
Even though the formula~\eqref{eq:COM_formula} is stated for
deterministic functions, it is easily extended to $\clZ_t$-measurable
functions which is how it is used above.  Therefore,
\begin{align*}
	\sJ_T(U) & = \var_0(Y_0) + \E\Big(\int_0^T (\Gamma Y_t)(X_t) +
                   |U_t + V_t(X_t)|^2 \ud t\Big)\\
& = \var_0(Y_0) + \tE\Big(\int_0^T \sigma_t(\Gamma Y_t) + \sigma_t(|U_t
  + V_t|^2) \ud t \Big)\\
& = \var_0(Y_0) + \tE\Big(\int_0^T \ell(Y_t,V_t,U_t;\sigma_t)\ud t\Big)
\end{align*}


\subsection{Proof of~\Thm{thm:optimal-solution}}\label{ss:pf-thm52}

Equation~\eqref{eq:Hamilton_eqns} is the Hamilton's equation for
optimal control of a BSDE~\cite[Thm.~4.4]{peng1993backward}. The
optimal control is obtained from the maximum principle:
\[
U_t = \mathop{\operatorname{argmax}}_{u\in\Re^m} \; \clH(Y_t,V_t,u,P_t;\sigma_t)
\]
Since $\clH$ is quadratic in the control input, the explicit
formula~\eqref{eq:opt-cont-soln} is obtained by evaluating the derivative and setting it to zero:
\[
\clH_u(Y_t,V_t,u,P_t;\sigma_t) = 2\sigma_t(\ones)u + 2\sigma_t(V_t) + P_t(h) = 0
\]

\subsection{Justification of the formula~\eqref{eq:P-t-ansatz}}\label{ssec:pf-ansatz}

For notational ease, we drop the superscript $^\opt$ and denote the
optimal control input simply as $U_t$.  In this proof, $\langle
\cdot,\cdot\rangle$ is used to denote the duality paring between
functions and measures (e.g., $\langle f, \mu\rangle = \mu(f)$).  

Let $f$ be an arbitrary test
function.  We show that
\[
\langle f, P_t\rangle = \big\langle 2f(Y_t-\pi_t(Y_t)), \sigma_t
\big\rangle,\quad 0< t \leq T
\]
This is known to be true at time $t=0$ because of the boundary
condition~\eqref{eq:Hamilton_eqns2-c}.  Therefore, the proof is
carried out by taking a derivative of both sides and showing these to
be identical.     

Using the It\^o-Wentzell formula for
measure valued processes~\cite[Thm.~1.1]{krylov2011ito}, 
\begin{align*}
\ud \big\langle 2f(Y_t&-\pi_t(Y_t)), \sigma_t 
\big\rangle \\ & = 2\big\langle \clA(fY_t)-f(\clA Y_t) -
                 \pi_t(Y_t)(\clA f),\sigma_t\big\rangle \ud t \\
 & \qquad\quad + (\langle 2f(U_t+V_t),\sigma_t \rangle+ \langle
  fh,P_t\rangle) \ud Z_t 
\end{align*}
where we have used $\ud
\big(\pi_t(Y_t)\big) = - U_t \ud Z_t$
(Prop.~\ref{prop:optimal-solution}).  From the Hamilton's
equation~\eqref{eq:Hamilton_eqns2-b}, upon explicitly evaluating the terms
\begin{align*}
	\ud \langle f, P_t\rangle &= \Big(\langle \clA f, P_t \rangle
                                    + \frac{\ud}{\ud \epsilon}\sigma_t\big(\Gamma (Y_t+\epsilon f)\big)\Big|_{\epsilon = 0} \Big)\ud t \\
	&\qquad\quad + (\langle fh, P_t \rangle + \langle 2f(U_t + V_t), \sigma_t \rangle)\ud Z_t
\end{align*}
where 
\[
\frac{\ud}{\ud \epsilon}\Gamma (Y_t+\epsilon f)\Big|_{\epsilon = 0} = 2\big(\clA(Y_tf)-Y_t(\clA f) - f(\clA Y_t)\big)
\] 
On comparing the terms, the two derivatives are the seen to be the
same where we use also the identity $\langle g, P_t\rangle = \big\langle 2g(Y_t-\pi_t(Y_t)), \sigma_t
\big\rangle$ for $g = \clA f$. 
%
%

\subsection{Proof of~\Thm{thm:martingale}}\label{ss:pf-thm53}

The proof uses the equation of the nonlinear filter and $\ud I_t := \ud
Z_t - \pi_t(h) \ud t$ is the innovation increment.  
We evaluate the derivative of $\clV_t(Y_t) = \pi_t(Y_t^2) - \big(\pi_t(Y_t)\big)^2$.
\begin{align*}
	\ud \pi_t(&Y_t^2)\\
	 &= \pi_t(\clA Y_t^2) \ud t + \big(\pi_t(hY_t^2)-\pi_t(h)\pi_t(Y_t^2)\big)\ud I_t \\
	&\quad+\pi_t\big(-2Y_t\big(\clA Y_t + h (U_t+V_t)\big) + |V_t|^2\big)\ud t \\
	&\quad+ 2\pi_t\big(Y_tV_t\big) \ud Z_t +2\big(\pi_t(h Y_tV_t)-\pi_t(h)\pi_t(Y_tV_t)\big)\ud t\\
	&=\pi_t\big(\Gamma Y_t\big) \ud t + \pi_t(|V_t|^2)\ud t - 2\pi_t(hY_t)U_t \ud t\\
	&\quad+\big(\pi_t(hY_t^2)-\pi_t(h)\pi_t(Y_t^2) + 2\pi_t(Y_tV_t)\big)\ud I_t
\end{align*}
Similarly, 
\begin{align*}
	\ud \pi_t(Y_t) &= \pi_t(\clA Y_t) \ud t \\
	&\quad+ \big(\pi_t(hY_t)-\pi_t(h)\pi_t(Y_t)\big)\big(\ud Z_t - \pi_t(h)\ud t\big) \\
	&\quad-\pi_t\big(\clA Y_t + h(U_t+V_t)\big)\ud t + \pi_t\big(V_t\big) \ud Z_t \\
	&\quad+\big(\pi_t(h V_t)-\pi_t(h) \pi_t(V_t)\big)\ud t \\
	&= \big(\pi_t(hY_t)-\pi_t(h)\pi_t(Y_t)+ \pi_t(V_t)\big)\ud Z_t \\
	&\quad- \big(U_t+\pi_t(hY_t)-\pi_t(h)\pi_t(Y_t)+ \pi_t(V_t)\big) \pi_t(h)\ud t \\
	&= U_t^\opt \ud Z_t - (U_t-U_t^\opt)\pi_t(h)\ud t
\end{align*}
where $U_t^\opt := -\pi_t(hY_t)+\pi_t(h)\pi_t(Y_t) - \pi_t(V_t)$. 
Therefore,
\begin{align*}
	\ud \big(\pi_t(Y_t)\big)^2 &= 2\pi_t(Y_t)U_t^\opt \ud Z_t  + |U_t^\opt|^2 \ud t\\
	&\quad - 2 \pi_t(Y_t)(U_t-U_t^\opt)\pi_t(h)\ud t
\end{align*}
Collecting terms, we have
\begin{align*}
	\ud M_t =& \pi_t\big(\Gamma Y_t\big) \ud t + \pi_t(|V_t|^2)\ud t - 2\pi_t(hY_t)U_t \ud t\\
	& +\big(\pi_t(hY_t^2)-\pi_t(h)\pi_t(Y_t^2) + 2\pi_t(Y_tV_t)\big)\ud I_t\\
	&-2\pi_t(Y_t)U_t^\opt \ud Z_t + 2 \pi_t(Y_t)\big(U_t-U_t^\opt\big)\pi_t(h)\ud t \\
	&- |U_t^\opt|^2 \ud t - \ell(Y_t,V_t,U_t;\pi_t) \ud t\\
	=&\big(\pi_t(hY_t^2)-\pi_t(h)\pi_t(Y_t^2) + 2\pi_t(Y_tV_t)\big)\ud I_t\\
	&-|U_t-U_t^\opt|^2 \ud t 
\end{align*}
Since $-|U_t - U_t^\opt|^2 \le 0$ and $I$ is a $\sP$-martingale, $M$ is a $\sP$-supermartingale, and it is a martingale if and only if $U_t = U_t^\opt$ for all $t$. 

\subsection{Formal derivation of the nonlinear filter}\label{ssec:pf-nonlinear-filter-derivation}


We begin with an ansatz
\begin{equation}\label{eq:pit_ansatz}
\ud \pi_t(f) = \alpha_t(f) \ud t + \beta_t(f) \ud Z_t
\end{equation}
where the goal is to obtain formulae for $\alpha_t$ and $\beta_t$.  Because we have 
an equation~\eqref{eq:estimator-with-optimal-ctrl} for $\pi_t(Y_t)$,
let us express $\ud(\pi_t(Y_t))$ in terms of the unknown $\alpha_t$
and $\beta_t$.  Using the SDE~\eqref{eq:pit_ansatz} for $\pi_t$ and the
BSDE~\eqref{eq:dual-optimal-control-b} for $Y_t$, apply the
It\^o-Wentzell formula to obtain
\begin{align*}
	\ud \big(\pi_t (Y_t)\big) &=
	( \alpha_t(Y_t) + \beta_t(V_t) - \pi_t(\clA Y_t + h^\tp (U_t+V_t)) )\ud t\\
	&\qquad  + (\beta_t(Y_t) + \pi_t(V_t))\ud Z_t \\
\end{align*}
Comparing with~\eqref{eq:estimator-with-optimal-ctrl}, 
\begin{align*}
\alpha_t(Y_t) + \beta_t(V_t) - \pi_t(\clA Y_t + h^\tp (U_t+V_t)) = 0 \\
\beta_t(Y_t) + \pi_t(V_t)  = \big(\pi_t(hY_t)-\pi_t(h)\pi_t(Y_t)\big)
  + \pi_t(V_t)
\end{align*}
for $ 0\le t \le T$, $\sP\text{-a.s.}$.  Because $F$ and therefore
$Y_t$ is arbitrary, the second of these equations suggests setting
\[
\beta_t(f) = \pi_t(hf) - \pi_t(h)\pi_t(f)
\]
using which the first equation is manipulated to show
\begin{align*}
	\alpha_t(Y_t) &= \pi_t(\clA Y_t) - \pi_t(h)\big(\pi_t(hY_t)-\pi_t(h)\pi_t(Y_t)+\pi_t(V_t)\big)\\
	&\qquad \quad + \pi_t(hV_t) -  \pi_t(hV_t)+\pi_t(h)\pi_t(V_t)\\
	&=\pi_t(\clA Y_t) - \pi_t(h)\big(\pi_t(hY_t)-\pi_t(h)\pi_t(Y_t)\big)
\end{align*}
which gives the following 
\[
\alpha_t(f) = \pi_t(\clA f) - \beta_t(f)\pi_t(h)
\]
Substituting the expressions for $\alpha_t$ and $\beta_t$ into the ansatz~\eqref{eq:pit_ansatz}
\begin{align*}
	&\ud \pi_t(f) 
	= \big(\pi_t(\clA f) - \beta_t(f)\pi_t(h)\big)\ud t + \beta_t(f)\ud Z_t\\
	&\;\;=\pi_t(\clA f) \ud t + \big(\pi_t(hf) - \pi_t(h)\pi_t(f)\big)(\ud Z_t - \pi_t(h)\ud t)
\end{align*}
This is the well known SDE of the nonlinear filter.

\subsection{Mitter-Newton duality for the linear Gaussian model} \label{apdx:Mitter-LG-case}

Consider~\eqref{eq:tildeX-MN-G} with the linear control
law~\eqref{eq:U-MN-G}.  Then $\tilde{X}_t$ is a Gaussian random
variable whose mean $\tilde{m}_t$ and variance $\tilde{\Sigma}_t$
evolve as follows:
\begin{subequations}
\begin{align}
	\frac{\ud \tilde{m}_t}{\ud t} &= A^\tp \tilde{m}_t + \sigma {u}_t\label{eq:tilde_mean}\\
	\frac{\ud \tilde{\Sigma}_t}{\ud t} &= (A^\tp+\sigma K_t) \tilde{\Sigma}_t +\tilde{\Sigma}_t(A^\tp+\sigma K_t)^\tp + \sigma\sigma^\tp \label{eq:tilde_var}
\end{align}
\end{subequations}
Note that the two equations are entirely un-coupled: $u_t$ affects
only the equation for $\tilde{m}_t$ and $K_t$ affects only the
equation for $\tilde{\Sigma}_t$.  We now turn to
explicitly computing the running cost.  For the linear Gaussian model
\[
(\clA^u h)(x) = H^\tp (A^\tp x + \sigma u)
\]
and the running cost becomes 
\begin{align*}
	l(x,u\,;z_t) & = \half|u|^2 + |H^\tp x|^2 + z_t H^\tp (A^\tp x + \sigma u)
\end{align*}
Because $\tilde{X}_t \sim N(\tilde{m}_t,\tilde{\Sigma}_t)$, 
\begin{align*}
	\E\big(l(\tilde{X}_t,u_t\,;z_t)\big) &= \half |{u}_t|^2 + \half\tr(K_t^\tp K_t\tilde{\Sigma}_t) + \half|H^\tp \tilde{m}_t|^2 \\
	&\quad + \half\tr(HH^\tp\tilde{\Sigma}_t)+ z_t H^\tp(A^\tp \tilde{m}_t+\sigma {u}_t)
\end{align*}
and because $\tilde{\mu}$ from $\mu$ are both Gaussian, the divergence 
\begin{align*}
	\E\Big(\log \frac{\ud \tilde{\mu}}{\ud \mu}(\tilde{X}_0)\Big) &= \half(m_0-\tilde{m}_0)^\tp \Sigma_0^{-1}(m_0-\tilde{m}_0)\\
	&\quad + \half\log\frac{\det(\tilde{\Sigma}_0)}{\det(\Sigma_0)} - \frac{d}{2} + \half \tr(\tilde{\Sigma}_0\Sigma_0^{-1})
\end{align*}
and because $h(\cdot)$ is linear the terminal condition term
\[
\E\big(z_T h(\tilde{X}_T)\big) = z_TH^\tp \tilde{m}_T
\]
Combining all of the above, upon a formal integration by parts, $\sJ(\tilde{\mu},U;z)$ is
expressed as sum of two un-coupled costs
\begin{align*}
	\sJ_1(\tilde{m}_0,{u}\,;z) &= \half(m_0-\tilde{m}_0)^\tp \Sigma_0^{-1}(m_0-\tilde{m}_0) \\
	&+ \int_0^T \half |u_t|^2 + \half|\dot{z}_t - H^\tp \tilde{m}_t|^2\ud t \\
	\sJ_2(\tilde{\Sigma}_0,K\,;z) &= \half\log(\det(\tilde{\Sigma}_0))+ \half \tr(\tilde{\Sigma}_0\Sigma_0^{-1}) \\
	&+ \int_0^T \half\tr(K_t^\tp K_t\tilde{\Sigma}_t) +  \half\tr(HH^\tp\tilde{\Sigma}_t)\ud t
\end{align*}
plus a few constant terms that are not affected by the decision
variables.  The first of these costs subject to the ODE
constraint~\eqref{eq:tilde_mean} for the mean $\tilde{m}_t$ is the
classical minimum energy duality.

\bibliographystyle{IEEEtran}
\bibliography{../bibfiles/_master_bib_jin.bib,../bibfiles/jin_papers.bib,../bibfiles/extrabib.bib,../bibfiles/estimator_controller.bib}

\begin{IEEEbiography}[{\includegraphics[width=1in,height=1.25in,clip,keepaspectratio]{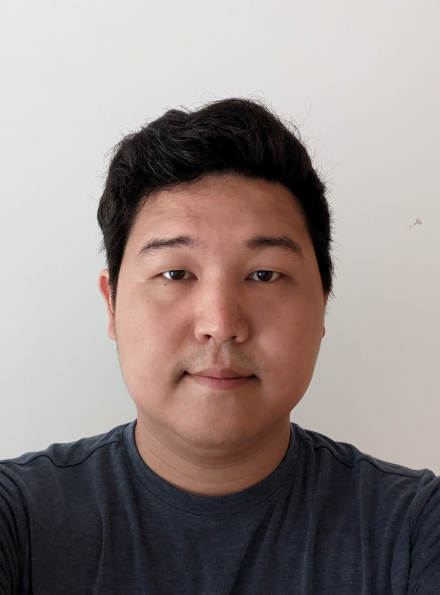}}]{Jin Won Kim} received the Ph.D. degree in Mechanical Engineering from University of Illinois at Urbana-Champaign, Urbana, IL, in 2022.
	He is now a postdocdoral research scientist in the Institute of Mathematics at the University of Potsdam.
	His current research interests are in nonlinear filtering and stochastic optimal control.
	He received the Best Student Paper Awards at the IEEE Conference on Decision and Control 2019.
\end{IEEEbiography}

\begin{IEEEbiography}[{\includegraphics[width=1in,height=1.25in,clip,keepaspectratio]{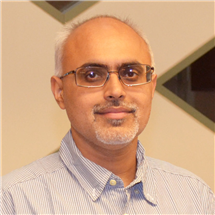}}]{Prashant G. Mehta} received the Ph.D. degree in Applied Mathematics from Cornell University, Ithaca, NY, in 2004.
	He is a Professor of Mechanical Science and Engineering at the University of Illinois at Urbana-Champaign.
	Prior to joining Illinois, he was a Research Engineer at the United Technologies Research Center (UTRC). His current research interests are in nonlinear filtering. He received the Outstanding Achievement Award at UTRC for his contributions to the modeling and control of combustion instabilities in jet-engines. His students received the Best Student Paper Awards at the IEEE Conference on Decision and Control 2007, 2009 and 2019, and were finalists for these awards in 2010 and 2012. In the past, he has served on the editorial boards of the ASME Journal of Dynamic Systems, Measurement, and Control and the Systems and Control Letters. He currently serves on the editorial board of the IEEE Transactions on Automatic Control. 	
\end{IEEEbiography}

\end{document}